\documentclass[12pt]{article}
\usepackage{amssymb}

\def\hybrid{\topmargin 0pt      \oddsidemargin 0pt
        \headheight 0pt \headsep 0pt
        \voffset=-0.5cm
        \hoffset=-0.25in
        \textwidth 6.75in
        \textheight 9.5in       
        \marginparwidth 0.0in
        \parskip 5pt plus 1pt   \jot = 1.5ex}
\catcode`\@=11
\def\marginnote#1{}

\newcount\hour
\newcount\minute
\newtoks\amorpm
\hour=\time\divide\hour by60 \minute=\time{\multiply\hour by60
\global\advance\minute by-\hour}
\edef\standardtime{{\ifnum\hour<12 \global\amorpm={am}%
        \else\global\amorpm={pm}\advance\hour by-12 \fi
        \ifnum\hour=0 \hour=12 \fi
        \number\hour:\ifnum\minute<10 0\fi\number\minute\the\amorpm}}
\edef\militarytime{\number\hour:\ifnum\minute<10 0\fi\number\minute}

\def\draftlabel#1{{\@bsphack\if@filesw {\let\thepage\relax
   \xdef\@gtempa{\write\@auxout{\string
      \newlabel{#1}{{\@currentlabel}{\thepage}}}}}\@gtempa
   \if@nobreak \ifvmode\nobreak\fi\fi\fi\@esphack}
        \gdef\@eqnlabel{#1}}
\def\@eqnlabel{}
\def\@vacuum{}
\def\draftmarginnote#1{\marginpar{\raggedright\scriptsize\tt#1}}
\def\draftlabel#1{{\@bsphack\if@filesw {\let\thepage\relax
   \xdef\@gtempa{\write\@auxout{\string
      \newlabel{#1}{{\@currentlabel}{\thepage}}}}}\@gtempa
   \if@nobreak \ifvmode\nobreak\fi\fi\fi\@esphack}
        \gdef\@eqnlabel{#1}}
\def\@eqnlabel{}
\def\@vacuum{}
\def\draftmarginnote#1{\marginpar{\raggedright\scriptsize\tt#1}}

\def\draft{\oddsidemargin -.5truein
        \def\@oddfoot{\sl preliminary draft \hfil
        \rm\thepage\hfil\sl\today\quad\militarytime}
        \let\@evenfoot\@oddfoot \overfullrule 3pt
        \let\label=\draftlabel
        \let\marginnote=\draftmarginnote
   \def\@eqnnum{(\theequation)\rlap{\kern\marginparsep\tt\@eqnlabel}%
\global\let\@eqnlabel\@vacuum}  }

\def\numberbysection{\@addtoreset{equation}{section}
        \def\theequation{\thesection.\arabic{equation}}}

\def\underline#1{\relax\ifmmode\@@underline#1\else
        $\@@underline{\hbox{#1}}$\relax\fi}

\def\titlepage{\@restonecolfalse\if@twocolumn\@restonecoltrue\onecolumn
     \else \newpage \fi \thispagestyle{empty}\c@page\z@
        \def\thefootnote{\fnsymbol{footnote}} }

\def\endtitlepage{\if@restonecol\twocolumn \else  \fi
        \def\thefootnote{\arabic{footnote}}
        \setcounter{footnote}{0}}  

 \hybrid

\newcounter{mo}

\newcommand{\vth}{\vartheta}

\newcommand{\Mat}{ {\rm Mat}(N,\mathbb C) }
\newcommand{\MatM}{ {\rm Mat}(M,\mathbb C) }

\newcommand{\mC}{\mathbb C}

\def\beq{\begin{equation}}
\def\eq{\end{equation}}
\def\p{\partial}

\begin{document}

\setcounter{page}{1}


\begin{flushright}
\end{flushright}
\vspace{0mm}

\begin{center}
\vspace{-10mm}
{\Large{GL(NM) quantum dynamical $R$-matrix based on solution}}\\
 \vspace{4mm}
 {\Large{of the associative Yang-Baxter equation}}
\\
\vspace{15mm} {\large {I. Sechin}
 \ \ \ \ \ \ \ \ \ \ \ {A. Zotov}
  }
 \vspace{10mm}

%
{\small{\rm
 Steklov Mathematical Institute of Russian Academy of Sciences,\\ Gubkina str. 8, Moscow,
119991,  Russia}}

\end{center}

\begin{center}\footnotesize{{\rm E-mails:}{\rm\
shnbuz@gmail.com,\ zotov@mi-ras.ru}}\end{center}
%
%

 \begin{abstract}
In this letter we construct ${\rm GL}_{NM}$-valued dynamical
$R$-matrix by means of unitary skew-symmetric  solution of the
associative Yang-Baxter equation in the fundamental representation
of ${\rm GL}_{N}$. In $N=1$ case the obtained answer reproduces the
${\rm GL}_{M}$-valued Felder's $R$-matrix, while in the $M=1$ case
it provides the
 ${\rm GL}_{N}$ $R$-matrix of vertex type including the
 Baxter-Belavin's
elliptic one and its degenerations.
 \end{abstract}



\vskip2cm

\paragraph{Yang-Baxter equations.} Consider a matrix-valued function $R_{12}^\hbar(z)\in\Mat^{\otimes
2}$, which solves the \underline{associative Yang-Baxter equation}
\cite{FK,Pol}:
  \beq\label{a01}
  R^\hbar_{12}(z_{12})
 R^{\eta}_{23}(z_{23})=R^{\eta}_{13}(z_{13})R_{12}^{\hbar-\eta}(z_{12})+
 R^{\eta-\hbar}_{23}(z_{23})R^\hbar_{13}(z_{13})\,,\quad
 z_{ab}=z_a-z_b\,.
  \eq
 Here, following notations of the Quantum Inverse Scattering Method
 \cite{Skl0}, an operator
  $R_{ab}^\hbar(z)$ in (\ref{a01})
   is considered as $\Mat^{\otimes 3}$-valued. It acts non-trivially in the $a$-th and $b$-th tensor components
   only. For example, $R_{13}^\hbar(z)$ is of the form
  \beq\label{a011}
  R^\hbar_{13}(z_{13})=\sum\limits_{i,j,k,l=1}^N R_{ijkl}(\hbar,z_{13})\, e_{ij}\otimes
  1_N\otimes e_{kl}\,,
  \eq
where the set $\{e_{ij}\}$ is the standard basis in $\Mat$, $1_N$ --
is the identity matrix in $\Mat$ and $R_{ijkl}(\hbar,z_{12})$ are
functions of complex variables $\hbar$ (the Planck constant) and $z$
(the spectral parameter).

    Let the solution of (\ref{a01}) satisfies also the properties of the
    \underline{skew-symmetry}
  \beq\label{a02}
  \begin{array}{c}
  \displaystyle{
 R^\hbar_{12}(z)=-R_{21}^{-\hbar}(-z)=-P_{12}R_{12}^{-\hbar}(-z)P_{12}\,,
 \qquad
 P_{12}=\sum\limits_{i,j=1}^N e_{ij}\otimes e_{ji}
 }
 \end{array}
 \eq
 and \underline{unitarity}
 \beq\label{a03}
   \begin{array}{c}
 \displaystyle{
R^\hbar_{12}(z) R^\hbar_{21}(-z) = (\wp(\hbar)-\wp(z))\,1_N\otimes
1_N\,,
 }
  \end{array}
  \eq
 where $\wp(x)$ -- is the Weierstrass
$\wp$-function. We assume that it is equal to $1/\sinh^2(x)$ or
$1/x^2$ for trigonometric (hyperbolic) or rational $R$-matrices
respectively. Notice that solution of (\ref{a01}) with the
properties (\ref{a02})-(\ref{a03}) is a true quantum $R$-matrix of
vertex type, i.e. it satisfies the \underline{quantum
(non-dynamical) Yang-Baxter equation}\footnote{The latter statement
is easily
 verified. See e.g. \cite{LOZ9}.}:
  \beq\label{a04}
R_{12}^\hbar(z_{12})R_{13}^\hbar(z_{13})R_{23}^\hbar(z_{23})=
R_{23}^\hbar(z_{23})R_{13}^\hbar(z_{13})R_{12}^\hbar(z_{12})\,.
  \eq
 Equation (\ref{a01}) can be view as matrix extension of
the \underline{genus one Fay trisecant identity}:
  \beq\label{a05}
  \begin{array}{c}
  \displaystyle{
\phi(\hbar,z_{12})\phi(\eta,z_{23})=\phi(\eta,z_{13})\phi(\hbar-\eta,z_{12})+\phi(\eta-\hbar,z_{23})\phi(\hbar,z_{13})\,,
 }
 \end{array}
 \eq
 which coincides with (\ref{a01}) in scalar ($N=1$) case.
  It plays a crucial role in the theory of classical and quantum integrable systems
  \cite{Krich1,Skl1,BFV}. Solution of (\ref{a05}) satisfying the (scalar versions
  of) properties (\ref{a02})-(\ref{a03}) is the \underline{Kronecker function}:
  \beq\label{a06}
  \begin{array}{c}
  \displaystyle{
\phi(\hbar,z)=\frac{\vth'(0)\vth(\hbar+z)}{\vth(\hbar)\vth(z)}\,,\quad
\vth(x)=\displaystyle{\sum _{k\in \mathbb Z}} \exp \left ( \pi
\imath \tau (k+\frac{1}{2})^2 +2\pi \imath
(x+\frac{1}{2})(k+\frac{1}{2})\right )\,,
 }
 \end{array}
 \eq
 where $\hbox{Im}(\tau)>0$. Its trigonometric and rational limits
 are given by $\coth(\hbar)+\coth(z)$ and $\hbar^{-1}+z^{-1}$
 respectively. Similarly, the elliptic solution of  (\ref{a01}) with properties
 (\ref{a02})-(\ref{a03}) is known \cite{Pol} to be given by the
 Baxter-Belavin's $R$-matrix \cite{Belavin}. The trigonometric
 solutions were classified in \cite{Sch}. They include the XXZ
 $R$-matrix, its 7-vertex deformation \cite{Chered} and their ${\rm
 GL}_N$ generalizations \cite{AHZ} (see a brief review in
 \cite{KrZ}). The rational solutions consist of the XXX
 $R$-matrix, its 11-vertex deformation \cite{Chered} and their
 ${\rm GL}_N$ generalizations \cite{Smirnov,LOZ8} -- deformations of
 the ${\rm GL}_N$
 Yang's  $R$-matrix $R_{12}^\hbar(z)=\hbar^{-1}1_N\otimes
 1_N+z^{-1}P_{12}$. Summarizing, we deal with the  $R$-matrices
 considered as matrix
 generalizations of the Kronecker function (including its
 trigonometric and rational versions).

To formulate the main result we also need the \underline{Felder's
dynamical ${\rm GL}_M$ $R$-matrix} \cite{Felder2}:
 \beq\label{a07}
 \begin{array}{c}
  \displaystyle{
 R^{\hbox{\tiny{F}}}_{12}(\hbar,z_1,z_2|\,q)=R^{\hbox{\tiny{F}}}_{12}(\hbar,z_1-z_2|\,q)=
 }
\\ \ \\
  \displaystyle{
 =\phi(\hbar,z_1-z_2)\sum\limits_{i=1}^M
 E_{ii}\otimes E_{ii}+\sum\limits_{i\neq j}^M
 E_{ij}\otimes E_{ji}\,
 \phi(z_1-z_2,q_{ij})
 +\sum\limits_{i\neq j}^M
 E_{ii}\otimes E_{jj}\, \phi(\hbar,-q_{ij})\,,
 }
 \end{array}
 \eq
 where $q_1,...,q_M$ -- are (free) dynamical parameters,
  \beq\label{a071}
  \begin{array}{c}
  \displaystyle{
 q_{ij}=q_i-q_j\,,
 }
 \end{array}
 \eq
  and the set $\{E_{ij}\}$ is the standard basis in $\MatM$.

 The $R$-matrix (\ref{a07}) is a solution of the \underline{quantum dynamical Yang-Baxter
 equation}:
  \beq\label{a08}
  \begin{array}{c}
  \displaystyle{
 R^\hbar_{12}(z_1,z_2|\,q)R^\hbar_{13}(z_1,z_3|\,q-\hbar^{(2)})R^\hbar_{23}(z_2,z_3|\,q)=\hspace{40mm}
}
\\ \ \\
  \displaystyle{
\hspace{40mm}
=R^\hbar_{23}(z_2,z_3|\,q-\hbar^{(1)})R^\hbar_{13}(z_1,z_3|\,q)R^\hbar_{12}(z_1,z_2|\,q-\hbar^{(3)})\,,
 }
 \end{array}
 \eq
where the shifts of the dynamical arguments $\{q_i\}$ are performed
as follows:
  \beq\label{a09}
  \begin{array}{c}
  \displaystyle{
R^\hbar_{12}(z_1,z_2|\,q+\hbar^{(3)})=P_3^\hbar\,
R^\hbar_{12}(z_1,z_2|\,q)\, P_3^{-\hbar} \,,\quad
P_3^\hbar=\sum\limits_{k=1}^M 1_M\otimes 1_M\otimes E_{kk}
\exp\Big(\hbar\frac{\p}{\p q_k}\Big)\,.
 }
 \end{array}
 \eq

\paragraph{Quantum dynamical ${\rm GL}_{NM}$ $R$-matrix.}
%
 %
 Consider the following ${\rm Mat}(NM,\mC)$-valued expression:
  \beq\label{a10}
  \begin{array}{c}
  \displaystyle{
  {\bf R}^\hbar_{1'2'12}(z,w)=\sum\limits_{i=1}^M\stackrel{1'}{E}_{ii}\otimes\stackrel{2'}{E}_{ii}
  \otimes\, R^\hbar_{12}(z-w)
+\sum\limits_{i\neq j}^M
\stackrel{1'}{E}_{ij}\otimes\stackrel{2'}{E}_{ji}\otimes
R_{12}^{\,q_{ij}}(z-w)+
 }
 \\
  \displaystyle{
+\sum\limits_{i\neq
j}^M\stackrel{1'}{E}_{ii}\otimes\stackrel{2'}{E}_{jj}\otimes
 \stackrel{1}{1}_{N}\otimes \stackrel{2}{1}_{N}\phi(\hbar,- q_{ij})\,,
 }
 \end{array}
 \eq
 where the ${\rm Mat}(NM,\mC)$ indices are represented
 in a way that the $\MatM$-valued tensor components
 are numbered by the primed numbers, and the $\Mat$-valued components are
those without primes (as previously). Put it differently, the
indices are arranged through
  ${\rm Mat}(NM,\mC)^{\otimes 2}\cong {\rm Mat}(M,\mC)^{\otimes 2}\otimes {\rm Mat}(N,\mC)^{\otimes 2}$.
  The order of tensor components is, in
fact, not important. It is chosen as in (\ref{a10}) just to
emphasize its similarity with the Felder's $R$-matrix (\ref{a07}).
The latter is reproduced from (\ref{a10}) in the $N=1$ case, when
the ${\rm GL}_N$ $R$-matrix entering (\ref{a10}) turns into the
Kronecker function (\ref{a06}).


The results of the paper are summarized in the following

\noindent {\bf Theorem}
 {\em
Let $R_{12}^\hbar(z)$ be some ${\rm GL}_N$ quantum non-dynamical
$R$-matrix satisfying the associative Yang-Baxter equation
(\ref{a01}) and the properties (\ref{a02})-(\ref{a03}). Then the
expression (\ref{a10}) is a quantum dynamical $R$-matrix, i.e. it
satisfies the quantum dynamical
 Yang-Baxter equation:
  \beq\label{a20}
  \begin{array}{c}
  \displaystyle{
    \mathbf{R}^\hbar_{1'2'12}(z_1, z_2 \mid q)
        \mathbf{R}^\hbar_{1'3'13}(z_1, z_3 \mid q - \hbar^{(2)})
            \mathbf{R}^\hbar_{2'3'23}(z_2, z_3 \mid q) =
            }
             \\ \ \\
   \displaystyle{
              =
    \mathbf{R}^\hbar_{2'3'23}(z_2, z_3 \mid q - \hbar^{(1)})
        \mathbf{R}^\hbar_{1'3'13}(z_1, z_3 \mid q)
            \mathbf{R}^\hbar_{1'2'12}(z_1, z_2 \mid q -
            \hbar^{(3)})\,,
  }
 \end{array}
 \eq
 where the shifts of arguments $\{q_i\}$ are performed similarly to
 (\ref{a09}):
  \beq\label{a201}
  \begin{array}{c}
  \displaystyle{
    \mathbf{R}^\hbar_{1'2'12}(z_1, z_2 \mid q +\hbar^{(3)})=
            {\bf P}_{3'}^\hbar\,\,\mathbf{R}^\hbar_{1'2'12}(z_1, z_2 \mid q)\, {\bf P}_{3'}^{-\hbar} \,,
            }
             \\ \ \\
   \displaystyle{
   {\bf P}_{3'}^\hbar=\sum\limits_{k=1}^M \stackrel{1'}{1}_{M}\otimes\stackrel{2'}{1}_{M}\otimes\stackrel{3'}{E}_{kk}\otimes
 \stackrel{1}{1}_{N}\otimes\stackrel{2}{1}_{N}\otimes\stackrel{3}{1}_{N}\exp\Big(\hbar\frac{\p}{\p q_k}\Big)\,.
  }
 \end{array}
 \eq
}
\noindent\underline{\em{Proof:}}\quad
 It is useful to write (\ref{a01}) as
 \begin{equation}\label{a22}
            R^\hbar_{ab}(z_{ab}) R^\eta_{bc}(z_{bc}) =
                R^{\eta - \hbar}_{bc}(z_{bc}) R^\hbar_{ac}(z_{ac}) +
                    R^\eta_{ac}(z_{ac}) R^{\hbar - \eta}_{ab}(z_{ab}),
        \end{equation}
where $a,b,c$ are distinct numbers from the set $\{1,2,3\}$.
Besides (\ref{a22}) and the properties (\ref{a02})-(\ref{a03}) the
proof of (\ref{a20}) uses the Yang-Baxter equation (\ref{a04}) for
the $\mathrm{GL}_N$ $R$-matrix and the following cubic relation:
        \begin{equation}\label{a23}
            R^\hbar_{ab}(z_{ab}) R^\eta_{ac}(z_{ac}) R^\hbar_{bc}(z_{bc}) -
                R^\eta_{bc}(z_{bc}) R^\hbar_{ac}(z_{ac}) R^\eta_{ab}(z_{ab}) =
                    R^{\hbar + \eta}_{ac}(z_{ac}) (\wp(\hbar) -
                    \wp(\eta))\,,
        \end{equation}
which is true under hypothesis of the theorem. If $\hbar=\eta$ it
reduces to (\ref{a04}). In the general case
 (\ref{a23}) leads (due to skew-symmetry of its r.h.s.) to
  \beq\label{x24}
  \begin{array}{c}
  \displaystyle{
 R^\eta_{ab} R^\hbar_{ac} R^\eta_{bc}+R^\hbar_{ab} R^\eta_{ac} R^\hbar_{bc}=R^\eta_{bc}
 R^\hbar_{ac} R^\eta_{ab}+R^\hbar_{bc}
 R^\eta_{ac} R^\hbar_{ab}\,,\quad
 R_{ab}^\hbar=R^\hbar_{ab}(z_a-z_b)\,,
 }
 \end{array}
 \eq
 known as the
Yang-Baxter equation with two Planck constants \cite{LOZ16}.
 The verification of (\ref{a20}) is a straightforward but cumbersome calculation. Consider, for example,
 the equation
 arising in the tensor component
$\stackrel{1'}{E}_{ij}\otimes\stackrel{2'}{E}_{kk}\otimes\stackrel{3'}{E}_{ji}$
 with $i \ne j \ne k \ne i$:
  \beq\label{a24}
  \begin{array}{c}
  \displaystyle{
    R^{q_{ik}}_{12}(z_{12}) R^{q_{kj}}_{13}(z_{13}) R^{q_{ik}}_{23}(z_{23}) +
        \phi(\hbar, q_{ik}) \phi(\hbar, q_{ki}) R^{q_{ij}}_{13}(z_{13}) =
            }
             \\ \ \\
   \displaystyle{
              =
            R^{q_{kj}}_{23}(z_{23}) R^{q_{ik}}_{13}(z_{13}) R^{q_{kj}}_{12}(z_{12}) +
                \phi(\hbar, q_{kj}) \phi(\hbar, q_{jk})
                R^{q_{ij}}_{13}(z_{13})\,.
   }
 \end{array}
 \eq
 To prove it one should use (\ref{a23}) written in the form
  \beq\label{a25}
  \begin{array}{c}
  \displaystyle{
    R^{q_{ik}}_{12}(z_{12}) R^{q_{kj}}_{13}(z_{13}) R^{q_{ik}}_{23}(z_{23}) -
        R^{q_{kj}}_{23}(z_{23}) R^{q_{ik}}_{13}(z_{13}) R^{q_{kj}}_{12}(z_{12}) =
        }
 \\ \ \\
  \displaystyle{
            =(\wp(q_{ik}) - \wp(q_{kj})) R^{q_{ik} + q_{kj}}_{13}(z_{13})=
            (\wp(q_{ik}) - \wp(q_{kj})) R^{q_{ij}}_{13}(z_{13})
    }
 \end{array}
 \eq
 and the well-known property of the
Kronecker function (scalar version of the unitarity)
  \beq\label{a26}
  \begin{array}{c}
  \displaystyle{
    \phi(\hbar, q_{ik}) \phi(\hbar, q_{ki}) = \wp(\hbar) -
    \wp(q_{ik})\,,\quad\quad
    \phi(\hbar, q_{kj}) \phi(\hbar, q_{jk}) = \wp(\hbar) -
    \wp(q_{kj})\,.
    }
 \end{array}
 \eq
 %
 The rest of the tensor components are verified similarly. $\blacksquare$

In the elliptic case, when $R_{12}^\hbar(z)$ is the Baxter-Belavin's
$R$-matrix, the result of the theorem is known \cite{LOSZ4}. Similar
results for the classical $r$-matrices were obtained previously by
P. Etingof and O. Schiffmann  \cite{ESch} and later in
\cite{LOSZ,ZS}, where the Hitchin type systems were described on the
Higgs bundles with non-trivial characteristic classes. Recently,
these type models appeared in the context of $R$-matrix valued Lax
pairs and quantum long-range spin chains \cite{GrSZ0,GrSZ}. In
\cite{LOSZ4} the answer (\ref{a10}) was verified explicitly in the
elliptic case without use of the associative Yang-Baxter equation.
In this respect the approach of this paper provides much simpler
proof. What is more important, the answer (\ref{a10}) is also valid
for all trigonometric and rational degenerations of the elliptic
$R$-matrix (satisfying the properties required in the Theorem). In
the light of results of \cite{GrSZ} the $R$-matrix (\ref{a10}) is
the one necessary for quantization of the (generalized) model of
interacting tops.

\paragraph{Classical ${\rm gl}_{NM}$ $r$-matrix.} As a by-product of the Theorem we
also get the classical dynamical Yang-Baxter equation for the
classical $r$-matrix of the generalized interacting tops
\cite{GrSZ}. Consider the classical limit of the ${\rm GL}_N$
$R$-matrix from the Theorem:
  \beq\label{a41}
  \begin{array}{c}
  \displaystyle{
  R_{12}^\hbar(z)=\hbar^{-1}1_N\otimes 1_N+r_{12}(z)+O(\hbar)\,.
 }
 \end{array}
 \eq
 The coefficient $r_{12}(z)$ is the classical $r$-matrix, and the
 quantum Yang-Baxter equation (\ref{a04}) reduces in the limit
 (\ref{a41}) to the \underline{classical (non-dynamical) Yang-Baxter
 equation}:
  \beq\label{a42}
  \begin{array}{c}
  \displaystyle{
  [r_{12}(z_{12}),r_{13}(z_{13})]+[r_{12}(z_{12}),r_{23}(z_{23})]+[r_{13}(z_{13}),r_{23}(z_{23})]=0\,.
 }
 \end{array}
 \eq
 Similarly, the classical dynamical $r$-matrix appears from
 (\ref{a07}) through (\ref{a41}). It satisfies the \underline{classical
 dynamical Yang-Baxter equation}:
  \beq\label{a43}
  \begin{array}{c}
  \displaystyle{
  [r_{12}(z_{12}),r_{13}(z_{13})]+[r_{12}(z_{12}),r_{23}(z_{23})]+[r_{13}(z_{13}),r_{23}(z_{23})]+
  }
  \\ \ \\
  \displaystyle{
  [\hat{\p}_1,r_{23}(z_{23})]-[\hat{\p}_2,r_{13}(z_{13})]+[\hat{\p}_3,r_{12}(z_{12})]=0\,,
 }
 \end{array}
 \eq
 which underlies the Poisson structure of the spin Calogero-Moser
 model \cite{BAB}. Here
  \beq\label{a44}
  \begin{array}{c}
  \displaystyle{
 \hat{\p}_3=\sum\limits_{k=1}^M 1_M\otimes 1_M\otimes
 E_{kk}\,\p_{q_k}\,,\quad\quad P_3^\hbar\stackrel{(\ref{a09})}{=}
 {1_M}^{\otimes 3}+\hbar\,\hat{\p}_3+O(\hbar^2)\,.
 }
 \end{array}
 \eq
 In the same way, starting from the quantum $R$-matrix (\ref{a10})
 one gets the classical $r$-matrix
  \beq\label{a30}
  \begin{array}{c}
  \displaystyle{
  {\bf r}_{1'2'12}(z)=\sum\limits_{i=1}^M\stackrel{1'}{E}_{ii}\otimes\stackrel{2'}{E}_{ii}
  \otimes\, r_{12}(z)
+\sum\limits_{i\neq j}^M
\stackrel{1'}{E}_{ij}\otimes\stackrel{2'}{E}_{ji}\otimes
R_{12}^{\,q_{ij}}(z)\,,
 }
 \end{array}
 \eq
and the classical dynamical Yang-Baxter equation follows from
(\ref{a20}):
  \beq\label{a31}
  \begin{array}{c}
  \displaystyle{
    [\mathbf{r}_{1'2'12}(z_{12}), \mathbf{r}_{1'3'13}(z_{13})] +
        [\mathbf{r}_{1'2'12}(z_{12}), \mathbf{r}_{2'3'23}(z_{23})] +
            [\mathbf{r}_{1'3'13}(z_{13}), \mathbf{r}_{2'3'23}(z_{23})] +
            }
             \\ \ \\
   \displaystyle{
            +
                [\hat\partial_{1'}, \mathbf{r}_{2'3'23}(z_{23})] -
                [\hat\partial_{2'}, \mathbf{r}_{1'3'13}(z_{13})] +
                [\hat\partial_{3'}, \mathbf{r}_{1'2'12}(z_{12})] = 0\,.
    }
 \end{array}
 \eq
 with
  \beq\label{a33}
  \begin{array}{c}
  \displaystyle{
 \hat{\p}_{3'}=\sum\limits_{k=1}^M \stackrel{1'}{1}_{M}\otimes\stackrel{2'}{1}_{M}\otimes\stackrel{3'}{E}_{kk}\otimes
 \stackrel{1}{1}_{N}\otimes\stackrel{2}{1}_{N}\otimes\stackrel{3}{1}_{N}\p_{q_k}\,,\quad
 {\bf P}_{3'}^\hbar\stackrel{(\ref{a201})}{=}
 {1_{MN}}^{\otimes 3}+\hbar\,\hat{\p}_{3'}+O(\hbar^2)\,.
 }
 \end{array}
 \eq
%


\paragraph{Acknowledgments.}
The second author is a Young Russian Mathematics award winner. The
work was performed at the Steklov Mathematical Institute of Russian
Academy of Sciences, Moscow. This work is supported by the Russian
Science Foundation under grant 19-11-00062.



\begin{small}
 
\end{small}

\end{document}